\documentclass[11pt,a4paper]{article}
\usepackage[T1]{fontenc}
\usepackage[utf8]{inputenc}
\usepackage{lmodern}
\usepackage[left=23mm,right=23mm,top=20mm,bottom=20mm]{geometry}
\usepackage{amsmath,amssymb,amsthm,mathtools}
\usepackage{enumitem,microtype,booktabs}
\makeatletter
\renewcommand\section{\@startsection{section}{1}{\z@}{-2.2ex plus -.6ex minus -.2ex}{1.1ex plus .2ex}{\normalfont\large\bfseries}}
\makeatother
\setlist{topsep=3pt,itemsep=1pt,parsep=0pt}
\usepackage[numbers,sort&compress]{natbib}
\usepackage{hyperref}
\hypersetup{hidelinks,pdftitle={Diagonal and branchwise criteria for universal Fourier series},pdfauthor={Sergo A. Episkoposian}}
\newtheoremstyle{compactplain}{4pt plus 2pt minus 1pt}{4pt plus 2pt minus 1pt}{\itshape}{}{\bfseries}{.}{.5em}{}
\theoremstyle{compactplain}
\newtheorem{theorem}{Theorem}[section]
\newtheorem{lemma}[theorem]{Lemma}
\newtheorem{proposition}[theorem]{Proposition}
\newtheorem{corollary}[theorem]{Corollary}
\theoremstyle{remark}\newtheorem{remark}[theorem]{Remark}
\newcommand{\Es}{\mathcal E_s}
\newcommand{\U}{\mathcal U_{\mathrm{str}}}
\newcommand{\Uf}{\mathcal U_{\mathrm{fin}}}
\newcommand{\M}{\mathcal M}

\numberwithin{equation}{section}
\title{Diagonal and branchwise criteria for\\universal Fourier series}
\author{Sergo A. Episkoposian\\[3pt]
\small National Polytechnic University of Armenia\\
\small 105 Teryan St., Yerevan 0009, Armenia\\
\small \texttt{sergo.episkoposyan@polytechnic.am}}
\date{}
\begin{document}
\maketitle
\begin{abstract}
We characterize diagonal profiles admitting universal real trigonometric Fourier series with factors in $\mathcal E_s=\bigcap_{p>s}\ell^p$ and $\ell^{s,q}$, $2\le s<\infty$, $1\le q\le\infty$. For nonincreasing profiles we obtain exact dyadic criteria and logarithmic endpoints. For independently regular odd/even sine/cosine branches we characterize simultaneous subseries, sign, and rearrangement universality for finite almost-everywhere targets. One function in every $L^p$, $1\le p<2$, supplies the coefficients. A separation example admits this realization although every compatible nonincreasing majorant with bounded relative cost belongs to $\ell^1$. Residue-class decompositions extend the criterion; a weakest-model condition permits extended-real rearrangement targets. A Walsh--Paley counterpart yields the same diagonal multiplier thresholds for ternary $\{0,\pm1\}$-universality, together with a Fourier--Walsh realization. Trigonometric realization uses Pogosyan's majorant theorem; the multiplier, factorization, projection-transfer, separation, and Walsh Fourier-origin arguments are proved here.
\end{abstract}
\setlength{\abovedisplayskip}{4pt plus 2pt minus 2pt}
\setlength{\belowdisplayskip}{4pt plus 2pt minus 2pt}
\setlength{\abovedisplayshortskip}{4pt plus 2pt}
\setlength{\belowdisplayshortskip}{6pt plus 2pt minus 2pt}

\section{Introduction}
Following Menshov's representation and Talalyan's rearrangement theory \cite{Menshov1941,Talalyan1960}, we prescribe the coordinate cost $d_n=A_nc_n$ of a universal trigonometric series, with fixed $A$ and $c$ in a specified sequence space. The coefficients must come from one integrable function; targets are reached by subseries, signs, or permutations.

For finite $s\ge2$ put
\begin{equation}\label{eq:spaces}
\Es=\bigcap_{p>s}\ell^p,\quad D_Ac=(A_nc_n)_{n\ge0},\quad
R_s=\begin{cases}\infty,&s=2,\\2s/(s-2),&s>2.\end{cases}
\end{equation}
Theorems~\ref{thm:monotone} and~\ref{thm:lorentz} give exact nonincreasing-profile criteria for $\Es$ and every $\ell^{s,q}$. Corollary~\ref{cor:lorentz-models} distinguishes their critical logarithmic thresholds and proves the finite-$q$ equality case by a double-logarithmic witness.

The passage from the global monotone criterion to independently regular nonmonotone branch profiles is a central point of the paper. Theorem~\ref{thm:branches} gives an exact criterion for independently regular odd/even sine/cosine branches without ordering their models, while Theorem~\ref{thm:separation} shows that this branchwise mechanism is not a disguised single-majorant argument: it yields realizations outside every compatible single-monotone-majorant construction with bounded relative cost. Corollary~\ref{cor:residues} extends the finite-target criterion to residue classes. Corollary~\ref{cor:extended} covers extended-real rearrangement targets when one model is dominated by the others; no such assertion is made for arbitrary independent models.

The trigonometric realization input \cite[Theorem~3]{Pogosyan1983} is isolated in Lemma~\ref{lem:pog}; it is the only external realization theorem used in the trigonometric sufficiency arguments. The multiplier identities, monotone factorization, Lorentz criteria, projection obstructions, branchwise transfer, residue decomposition, and separation theorem are proved independently of its proof. The contribution claimed here is the resulting diagonal and branchwise classification, its sharp endpoints, and the separation mechanism, not the underlying classical multiplier inequalities themselves. Related coefficient-decay results \cite{Korner1989,GalstyanOvsepyan1996} concern different universality, null-series, or correction formulations and are not used as inputs. Lorentz diagonal inequalities are classical \cite{Carando2009}. Abstract universal series \cite{Bayart2008,Koumandos2010} likewise impose coefficient restrictions in a different approximation framework; the present conclusions require almost-everywhere convergence of selected, sign-modified, or rearranged real Fourier series in a prescribed diagonal image. Pointwise variants \cite{Edge1970,Shkarin2009}, sign universality \cite{Grigoryan2020Notes,GrigoryanGaloyan2021,Avetisyan2023}, weighted representation \cite{GrigorianEpiskoposian2001,EpiskoposianJFA2006,EpiskoposianIJMMS2006}, and prescribed Walsh coefficients \cite{GevorkyanNavasardyan1999} use different convergence or coefficient requirements. In particular, density of sign-modified partial sums in $L^p$, $p<1$, supplies neither the prescribed majorant nor the almost-everywhere limits required here.

Sections~\ref{sec:prelim}--\ref{sec:mult} give reductions, Sections~\ref{sec:monotone}--\ref{sec:examples} (including~\ref{sec:lorentz}) diagonal criteria, Sections~\ref{sec:obstructions}--\ref{sec:separation} (especially~\ref{sec:branches}--\ref{sec:extended}) the branch theory, and Section~\ref{sec:walsh} the Walsh analogue.

\section{Target classes, sequence facts, and realization}\label{sec:prelim}
On $[0,2\pi]$ use the orthonormal real system
\begin{equation}\label{eq:system}
\Phi_0=(2\pi)^{-1/2},\qquad \Phi_{2k-1}=\pi^{-1/2}\cos kx,\qquad
\Phi_{2k}=\pi^{-1/2}\sin kx\quad(k\ge1).
\end{equation}
Let $\mathcal S$ consist of real measurable functions finite almost everywhere, and let $\mathcal S^\infty$ consist of extended-real measurable functions. All almost-everywhere limits refer to Lebesgue measure. A series $\sum d_n\Phi_n$ is $\delta$-universal if every $F\in\mathcal S$ is an almost-everywhere limit of a subseries with selectors $\delta_n\in\{0,1\}$. It is $\gamma$-universal if the same holds after choices of signs $\gamma_n\in\{-1,1\}$. In both modes the order is the displayed order in \eqref{eq:system}. It is $\sigma_{\rm fin}$-universal if every $F\in\mathcal S$ is obtained by a permutation of all indices, and $\sigma$-universal if this holds for every $F\in\mathcal S^\infty$. Convergence is required along all partial sums, not only a subsequence.

Write $\Uf$ for Fourier coefficient sequences having simultaneously the $\delta$, $\gamma$, and $\sigma_{\rm fin}$ properties, and $\U$ for those having the $\delta$, $\gamma$, and $\sigma$ properties. Thus
\begin{equation}\label{eq:classes}
\U\subset\Uf.
\end{equation}
In each definition one coefficient sequence is fixed before the target is chosen. For $r>0$, membership in $\ell^r$ means $\sum|c_n|^r<\infty$; no normed-space assertion is made for $r<1$. For $c\in c_0$, let $c^*$ list the nonzero absolute values in nonincreasing order, appending zeros when the support is finite. Define, for $1\le q<\infty$,
\begin{equation}\label{eq:lorentz-definition}
\|c\|_{s,q}^q=\sum_{n\ge0}(n+1)^{q/s-1}(c_n^*)^q,\qquad
\ell^{s,q}=\{c\in c_0:\|c\|_{s,q}<\infty\},
\end{equation}
and
\begin{equation}\label{eq:weak}
\|c\|_{s,\infty}=\sup_{n\ge0}(n+1)^{1/s}c_n^*,\qquad
\ell^{s,\infty}=\{c\in c_0:\|c\|_{s,\infty}<\infty\}.
\end{equation}
These expressions are used as quasi-norms; no isometric norm assertion is needed.
A linear space $X\subset c_0$ is solid if $b\in X$ and $|c_n|\le|b_n|$ imply $c\in X$. Rearrangement stability below means explicitly $c\in X\Rightarrow c^*\in X$, including sequences with infinitely many zeros.

\begin{lemma}[Sequence-space facts]\label{lem:spaces}
For $0<p<r<\infty$, $\ell^p\subset\ell^r$ and $\|c\|_r\le\|c\|_p$. The spaces $\Es$ and $\ell^{s,q}$, $1\le q\le\infty$, are solid linear subspaces of $c_0$, stable under decreasing rearrangement, subsequence restriction, finite assembly, and uniformly bounded repetitions. Also $\ell^{s,q}\subset\ell^{s,\infty}\subset\Es$.
\end{lemma}
\begin{proof}
Normalization by $\|c\|_p$ gives
\begin{equation}\label{eq:sequence-inclusion}
\sum|c_n/\|c\|_p|^r\le\sum|c_n/\|c\|_p|^p=1.
\end{equation}
The zero case is immediate. A null sequence has only finitely many coordinates above any positive threshold; successive attained maxima therefore list every nonzero absolute value. This justifies rearrangement, including infinitely many zeros. Solidity and restriction follow from coordinatewise comparison of order statistics. For $\Es$, \eqref{eq:sequence-inclusion} and the power sums prove all stability statements.

For an assembly of $L$ sequences, at most $L(k-1)$ coordinates exceed $\max_i(u^{(i)})^*_{k-1}$; hence $v^*_{Lk-1}\le\max_i(u^{(i)})^*_{k-1}$, $k\ge1$. Splitting the Lorentz sum into the $L$ residue classes of indices and using $(Lk)^{q/s-1}\asymp_L k^{q/s-1}$ gives, for $q<\infty$, $\|v\|_{s,q}^q\le C_{s,q,L}\sum_i\|u^{(i)}\|_{s,q}^q$, while $\|v\|_{s,\infty}\le C_{s,L}\max_i\|u^{(i)}\|_{s,\infty}$. The same counting estimate applies when every value is repeated at most $L$ times, proving bounded-repetition stability. Also $(c+d)^*_{2k}\le c_k^*+d_k^*$; grouping consecutive indices and using $(a+b)^q\le2^{q-1}(a^q+b^q)$ proves linearity. The finitely many initial coordinates are harmless, and zero insertion preserves $c^*$. For the weak space the equivalent distribution bound is
\begin{equation}\label{eq:distribution}
N_c(t):=\#\{|c_n|>t\}\le(\|c\|_{s,\infty}/t)^s,\quad
N_{c+d}(t)\le N_c(t/2)+N_d(t/2).
\end{equation}
To verify the equivalence, take $0<t<c_n^*$ and then $t\uparrow c_n^*$; conversely use the last coordinate above $t$. Thus \eqref{eq:distribution} proves weak-space linearity directly. For finite $q$, summing over $(n+1)/2\le k+1\le n+1$ gives $\|c\|_{s,q}^q\ge C_{s,q}(n+1)^{q/s}(c_n^*)^q$, since each $c_k^*\ge c_n^*$ and the weights on this interval sum to at least $C_{s,q}(n+1)^{q/s}$ (handle $n=0$ directly). Finally
\begin{equation}\label{eq:weak-inclusion}
\sum|c_n|^p\le\|c\|_{s,\infty}^p\sum(n+1)^{-p/s}<\infty\quad(p>s).
\end{equation}
This proves both inclusions using \eqref{eq:weak-inclusion}.
\end{proof}

The following realization statement is the precise consequence of \cite[Theorem~3]{Pogosyan1983} used throughout.
\begin{lemma}[External majorant realization input]\label{lem:pog}
If $\rho_n>0$, $\rho_n\downarrow0$, and $\sum\rho_n^2=\infty$, there is a real $f\in\bigcap_{1\le p<2}L^p[0,2\pi]$ whose coefficients
\begin{equation}\label{eq:pog-coefficients}
d_n=\int_0^{2\pi}f(x)\Phi_n(x)\,dx,\qquad |d_n|<\rho_n,
\end{equation}
belong to $\U$.
\end{lemma}
\begin{proof}[Derivation from the published theorem]
Apply \cite[Theorem~3, p.~211; Russian original, p.~185]{Pogosyan1983} to $\rho/2$. Its real system is \eqref{eq:system}, and its definitions include $\mathcal S^\infty$ in the rearrangement mode. It gives one function in all the stated $L^p$ spaces with simultaneous universality and $|d_n|\le\rho_n/2<\rho_n$, hence \eqref{eq:pog-coefficients}. The source is a research announcement; its realization theorem is used as a published external result.
\end{proof}

\begin{remark}[Role of the external realization theorem]\label{rem:pog-role}
Lemma~\ref{lem:pog} is used only to pass from a positive nonincreasing nonsquare-summable majorant to actual trigonometric Fourier coefficients with the stated universality properties. No proof of Pogosyan's theorem is claimed here. All sequence-space multiplier computations and the structural projection, branch, residue, and separation arguments below are independent of its proof. A self-contained replacement would require an additional majorant-preserving block construction controlling both selected blocks and all intermediate partial sums; such a construction is not asserted in this paper.
\end{remark}

\begin{remark}[Hausdorff--Young baseline]\label{rem:HY}
Hausdorff--Young \cite[Vol.~II, Ch.~XII]{Zygmund2002} gives $d\in E_2\subset\Es$ for $f\in\bigcap_{1<p<2}L^p$: interpolate $L^1\to\ell^\infty$, bounded by $\pi^{-1/2}$, and $L^2\to\ell^2$ to obtain $L^p\to\ell^{p'}$. Thus profiles bounded below already admit $\Es$ factors by Lemma~\ref{lem:pog}, with $\rho_n=(n+1)^{-1/2}$. Profiles tending to zero require additional division control.
\end{remark}

\begin{proposition}[Square-summability obstruction]\label{prop:l2}
If $d\in\ell^2$, its series has none of the four universality properties above.
\end{proposition}
\begin{proof}
For any finite set $I$, orthonormality in \eqref{eq:system} gives
\begin{equation}\label{eq:orthogonal}
\left\|\sum_{n\in I}e_nd_n\Phi_n\right\|_2^2=\sum_{n\in I}|e_nd_n|^2
\le\sum_{n\in I}|d_n|^2\quad(|e_n|\le1).
\end{equation}
Completeness gives $L^2$ convergence for every selection or sign choice. A permutation has the same $L^2$ sum: after it includes an initial segment, \eqref{eq:orthogonal} bounds its error by the square-summable tail. An $L^2$-convergent sequence has an almost-everywhere convergent subsequence: choose squared errors at most $2^{-3j}$; Chebyshev at $2^{-j}$ bounds the exceptional measures by $2^{-j}$, whose tail unions tend to zero. Thus any almost-everywhere limit equals an $L^2$ function. The finite target $x^{-1/2}$ on $(0,2\pi]$ is not in $L^2$ and excludes all modes.
\end{proof}

\section{Constructive multiplier and solid-space reductions}\label{sec:mult}
For sequence spaces $X,Y$, write
\begin{equation}\label{eq:mult-definition}
\M(X,Y)=\{A:Ac\in Y\text{ for every }c\in X\}.
\end{equation}

\begin{lemma}[Finite diagonal test]\label{lem:holder}
For $q>2$, put $r=2q/(q-2)$. Then
\begin{equation}\label{eq:holder}
\|Ac\|_2\le\|A\|_r\|c\|_q\quad(A\in\ell^r),
\end{equation}
and a bound $\|Ac\|_2\le C\|c\|_q$ on all finitely supported $c$ implies $A\in\ell^r$.
\end{lemma}
\begin{proof}
The identity $2/r+2/q=1$ gives \eqref{eq:holder} by H\"older. For a finite set $I$ set $v_n=|A_n|^{2/(q-2)}\mathbf1_I(n)$ and $S_I=\sum_{n\in I}|A_n|^r$. Then
\begin{equation}\label{eq:finite-test}
\|Av\|_2=S_I^{1/2},\qquad\|v\|_q=S_I^{1/q},\qquad
S_I^{1/r}\le C\quad(S_I>0).
\end{equation}
The bound is automatic when $S_I=0$. Increasing $I$ to initial segments in \eqref{eq:finite-test} gives $\sum|A_n|^r\le C^r$.
\end{proof}

\begin{theorem}[Exact multiplier threshold]\label{thm:mult}
For $2\le s<\infty$,
\begin{equation}\label{eq:multiplier}
\M(\Es,\ell^2)=\bigcup_{0<r<R_s}\ell^r.
\end{equation}
\end{theorem}
\begin{proof}
If $A\in\ell^{r_0}$ with $r_0<R_s$, choose $\max(2,r_0)<r<R_s$. Lemma~\ref{lem:spaces} gives $A\in\ell^r$. Defining $q>2$ by $1/q=1/2-1/r$ gives $q>s$, by \eqref{eq:spaces}. Hence \eqref{eq:holder} proves $Ac\in\ell^2$ for every $c\in\Es$.

For the converse we construct a witness, avoiding a continuity assumption. Suppose $A$ is outside the union in \eqref{eq:multiplier}. Set $q_j=s+1/j$ and $r_j=2q_j/(q_j-2)<R_s$. Every tail of $\sum|A_n|^{r_j}$ diverges. Choose successive disjoint finite intervals $I_j$ beyond the preceding interval, with
\begin{equation}\label{eq:gliding-blocks}
S_j:=\sum_{n\in I_j}|A_n|^{r_j}\ge2^{jr_j},\qquad
c_n^{(j)}=S_j^{-1/2}|A_n|^{2/(q_j-2)}\mathbf1_{I_j}(n).
\end{equation}
The finite computation of Lemma~\ref{lem:holder}, namely \eqref{eq:finite-test}, applied to \eqref{eq:gliding-blocks}, gives
\begin{equation}\label{eq:gliding-norms}
\|Ac^{(j)}\|_2=1,\qquad\|c^{(j)}\|_{q_j}=S_j^{-1/r_j}\le2^{-j}.
\end{equation}
Define $c$ by these values on the disjoint intervals and by zero elsewhere. For fixed $m$, $q_j\le q_m$ when $j\ge m$, so Lemma~\ref{lem:spaces} and \eqref{eq:gliding-norms} give
\begin{equation}\label{eq:gliding-conclusion}
\sum|c_n|^{q_m}=\sum_j\|c^{(j)}\|_{q_m}^{q_m}<\infty,
\qquad \sum|A_nc_n|^2=\sum_j1=\infty.
\end{equation}
The early blocks in \eqref{eq:gliding-conclusion} are finite. For every $q>s$ some $q_m<q$; Lemma~\ref{lem:spaces} therefore gives $c\in\Es$ but $Ac\notin\ell^2$, proving \eqref{eq:multiplier} by \eqref{eq:mult-definition}.
\end{proof}

\begin{lemma}[Rearrangement witness]\label{lem:rearrangement}
If $A_n\ge0$ is nonincreasing and $c\in c_0$, then
\begin{equation}\label{eq:rearrangement}
\sum_{n=0}^NA_n^2|c_n|^2\le\sum_{n=0}^NA_n^2(c_n^*)^2\quad(N\ge0).
\end{equation}
In particular $Ac\notin\ell^2$ implies $Ac^*\notin\ell^2$.
\end{lemma}
\begin{proof}
Sort the finite list $|c_0|,\ldots,|c_N|$ into $v_0\ge\cdots\ge v_N$. Swapping out-of-order values at $i<j$ increases the weighted sum by $(A_i^2-A_j^2)(|c_j|^2-|c_i|^2)\ge0$. Sorting therefore increases it. Each $v_n\le c_n^*$, since the finite list is a subset of the full sequence. This proves \eqref{eq:rearrangement}; increasing $N$ proves its last consequence, also when a sum is infinite.
\end{proof}

\begin{theorem}[Solid-space criterion]\label{thm:solid}
Let $X\subset c_0$ be a solid linear space satisfying $c\in X\Rightarrow c^*\in X$. For positive nonincreasing $A$,
\begin{equation}\label{eq:solid}
D_A(X)\cap\U\ne\varnothing\quad\Longleftrightarrow\quad A\notin\M(X,\ell^2).
\end{equation}
This criterion also characterizes existence in each mode separately, including $\sigma_{\rm fin}$. Sufficiency produces Fourier coefficients of one real function in every $L^p$, $1\le p<2$.
\end{theorem}
\begin{proof}
Necessity follows from \eqref{eq:mult-definition} and Proposition~\ref{prop:l2}. Otherwise choose $c\in X$ with $Ac\notin\ell^2$ and put $b=c^*$. Lemma~\ref{lem:rearrangement} gives $Ab\notin\ell^2$. There are infinitely many nonzero coordinates, since finite support would give $Ac\in\ell^2$. Thus $b_n>0$, $b\in X$, and $b_n\downarrow0$. The product $\rho=Ab$ is positive, nonincreasing, outside $\ell^2$, and tends to zero because $A_n\le A_0$. Lemma~\ref{lem:pog} gives $d$ with $|d_n|<A_nb_n$. Then
\begin{equation}\label{eq:quotient}
\widetilde c_n=d_n/A_n,\qquad |\widetilde c_n|\le b_n
\end{equation}
belongs to $X$ by solidity. Equation~\eqref{eq:quotient} gives $D_A\widetilde c=d$, proving \eqref{eq:solid} and all additional assertions.
\end{proof}

\section{Monotone factorization and diagonal thresholds}\label{sec:monotone}
\begin{lemma}[Divergent weighted tails]\label{lem:weighted}
Let $A_n>0$ be nonincreasing and outside $\bigcup_{0<r<R_s}\ell^r$. For every integer $j>s$,
\begin{equation}\label{eq:weighted}
\sum_{n\ge1}A_n^2n^{-2/s+2/j}=\infty.
\end{equation}
\end{lemma}
\begin{proof}
Suppose the sum in \eqref{eq:weighted} were finite, with value $S$. On $2^{m-1}\le n<2^m$, monotonicity and the negative exponent give
\begin{equation}\label{eq:dyadic}
S\ge\frac12\,2^{m(1-2/s+2/j)}A_{2^m}^2,
\quad A_n\le Cn^{-\alpha_j},\quad
\alpha_j=\frac12-\frac1s+\frac1j>0.
\end{equation}
The bound in \eqref{eq:dyadic} follows at $2^m$, then on $[2^m,2^{m+1})$ by $A_n\le A_{2^m}$ and $2^{-m}<2/n$. If $s=2$, it gives $A\in\ell^r$ for every $r>j$. If $s>2$, $\alpha_j>1/R_s$, so choosing $1/\alpha_j<r<R_s$ gives $A\in\ell^r$. Both contradict the hypothesis; constants may depend on $A,s,j$.
\end{proof}

\begin{theorem}[Explicit monotone factorization]\label{thm:factor}
For positive nonincreasing $A$, exclusion from $\bigcup_{0<r<R_s}\ell^r$ is equivalent to the existence of positive nonincreasing $b\in\Es$ with
\begin{equation}\label{eq:factor}
b_n\downarrow0,\qquad A_nb_n\downarrow0,\qquad Ab\notin\ell^2.
\end{equation}
The factor can additionally be chosen outside $\ell^s$.
\end{theorem}
\begin{proof}
The existence in \eqref{eq:factor} excludes the multiplier union by Theorem~\ref{thm:mult}. Conversely fix an integer $J>s$; the block index $j$ below runs through all integers $j\ge J$. Lemma~\ref{lem:weighted} and divergence of each tail allow successive integers $2\le N_J<N_{J+1}<\cdots$ such that
\begin{equation}\label{eq:factor-blocks}
\sum_{n=N_j}^{N_{j+1}-1}A_n^2n^{-2/s+2/j}\ge1,\qquad
b_n=n^{-1/s+1/j}\quad(N_j\le n<N_{j+1}).
\end{equation}
Set $b_n=b_{N_J}$ for $0\le n<N_J$. Inside blocks $b$ decreases, and at $N=N_{j+1}>1$,
\begin{equation}\label{eq:factor-boundary}
N^{-1/s+1/(j+1)}\le N^{-1/s+1/j}<(N-1)^{-1/s+1/j}.
\end{equation}
Thus \eqref{eq:factor-boundary} proves global monotonicity. Also $b_n\le n^{-1/s+1/J}\to0$. Given $q>s$, let $\eta=(q/s-1)/2>0$ and choose $j_0$ with $q/j\le\eta$ for $j\ge j_0$. Equation~\eqref{eq:factor-blocks} then gives
\begin{equation}\label{eq:factor-summability}
b_n^q\le n^{-1-\eta}\quad(n\ge N_{j_0}),\qquad
b_n^s\ge n^{-1}\quad(n\ge N_J).
\end{equation}
The first bound in \eqref{eq:factor-summability}, with finite earlier blocks, proves $b\in\Es$; the second proves $b\notin\ell^s$. Each block in \eqref{eq:factor-blocks} contributes at least one to the energy; $Ab$ decreases and $A_nb_n\le A_0b_n\to0$, proving \eqref{eq:factor}.
\end{proof}

\begin{theorem}[Monotone diagonal criterion]\label{thm:monotone}
For $2\le s<\infty$ and positive nonincreasing $A$,
\begin{equation}\label{eq:monotone}
D_A(\Es)\cap\U\ne\varnothing
\quad\Longleftrightarrow\quad
A\notin\bigcup_{0<r<R_s}\ell^r.
\end{equation}
The condition also characterizes existence in any one mode separately. The sufficient direction gives Fourier coefficients of one real function in every $L^p$, $1\le p<2$.
\end{theorem}
\begin{proof}
Theorems~\ref{thm:solid} and~\ref{thm:mult} give \eqref{eq:monotone}, including the individual modes and Fourier origin. Explicitly, Theorem~\ref{thm:factor} supplies \eqref{eq:factor}; Lemma~\ref{lem:pog} applied to $Ab$, followed by \eqref{eq:quotient}, gives the factor in $\Es$.
\end{proof}

\begin{theorem}[Weak-space diagonal criterion]\label{thm:weak}
For positive nonincreasing $A$,
\begin{equation}\label{eq:weak-criterion}
D_A(\ell^{s,\infty})\cap\U\ne\varnothing
\quad\Longleftrightarrow\quad
\sum_{n\ge0}A_n^2(n+1)^{-2/s}=\infty.
\end{equation}
The assertions about individual modes and Fourier origin in Theorem~\ref{thm:monotone} hold here as well.
\end{theorem}
\begin{proof}
If the sum is finite, Lemma~\ref{lem:rearrangement} and \eqref{eq:weak} give
\begin{equation}\label{eq:weak-bound}
\sum A_n^2|c_n|^2\le\sum A_n^2(c_n^*)^2
\le\|c\|_{s,\infty}^2\sum A_n^2(n+1)^{-2/s}<\infty.
\end{equation}
Proposition~\ref{prop:l2} excludes universality in \eqref{eq:weak-bound}. Otherwise $b_n=(n+1)^{-1/s}\in\ell^{s,\infty}$ and $Ab\notin\ell^2$. Theorem~\ref{thm:solid} gives \eqref{eq:weak-criterion}.
\end{proof}

\begin{corollary}[Comparable profiles]\label{cor:comparison}
Let $X$ be solid, $A_n\ne0$, $B_n>0$, and $mB_n\le|A_n|\le MB_n$ with $0<m\le M<\infty$. Then
\begin{equation}\label{eq:comparable-images}
D_A(X)=D_B(X).
\end{equation}
For $X=\Es$ or $\ell^{s,q}$ and nonincreasing $B$, the corresponding criterion in \eqref{eq:monotone} or \eqref{eq:lorentz-criterion} applies to $A$ with $B$ in its place.
\end{corollary}
\begin{proof}
If $d=Bv$, set $c_n=B_nv_n/A_n$; then $|c_n|\le|v_n|/m$, so $c\in X$ and $Ac=d$. Conversely $d=Ac$ can be written $Bu$ with $u_n=A_nc_n/B_n$ and $|u_n|\le M|c_n|$. Solidity and linearity prove both inclusions in \eqref{eq:comparable-images}. Theorems~\ref{thm:monotone} and~\ref{thm:lorentz} give the stated consequences.
\end{proof}

\section{Exact Lorentz criteria}\label{sec:lorentz}
For $A_n>0$ nonincreasing, define
\begin{equation}\label{eq:lorentz-weight}
w_j(A)=2^{j(1-2/s)}A_{2^j-1}^2\quad(j\ge0),\qquad
Y_q=\begin{cases}\ell^\infty,&1\le q\le2,\\
\ell^{q/(q-2)},&2<q<\infty,\\ \ell^1,&q=\infty.\end{cases}
\end{equation}
The shift $2^j-1$ in \eqref{eq:lorentz-weight} keeps the zero-based indexing of \eqref{eq:system}.

\begin{theorem}[Lorentz diagonal criterion]\label{thm:lorentz}
For $2\le s<\infty$, $1\le q\le\infty$, and positive nonincreasing $A$,
\begin{equation}\label{eq:lorentz-criterion}
A\in\M(\ell^{s,q},\ell^2)\iff w(A)\in Y_q,\qquad
D_A(\ell^{s,q})\cap\U\ne\varnothing\iff w(A)\notin Y_q.
\end{equation}
The second criterion holds for each universality mode separately. Sufficiency gives Fourier coefficients of one real function in every $L^p$, $1\le p<2$.
\end{theorem}
\begin{proof}
Use one-based coordinates $a_k=c_{k-1}^*$ and $\mathcal A_k=A_{k-1}$. Put $x_j=2^{j/s}a_{2^j}$. On $2^j\le k<2^{j+1}$, $a_{2^{j+1}}\le a_k\le a_{2^j}$, and the sum of $k^{q/s-1}$ is comparable to $2^{jq/s}$. Thus
\begin{equation}\label{eq:lorentz-dyadic}
\|c\|_{s,q}^q\asymp\sum_jx_j^q\ (q<\infty),\quad
\|c\|_{s,\infty}\asymp\sup_jx_j,\quad
\sum_n A_n^2|c_n|^2\le\sum_jw_j(A)x_j^2.
\end{equation}
For the lower norm bound the block endpoint gives $x_{j+1}$; the first coordinate supplies $x_0$. The energy bound uses Lemma~\ref{lem:rearrangement}. In \eqref{eq:lorentz-dyadic}, $\ell^q\subset\ell^2$ for $q\le2$, H\"older with exponents $q/2,q/(q-2)$ for $q>2$, and the supremum bound for $q=\infty$ show sufficiency of $w\in Y_q$ for the multiplier assertion.

For necessity the monotonicity constraint on $a$ must be retained. Given any nonnegative $x\in\ell^q$ (bounded if $q=\infty$), set
\begin{equation}\label{eq:lorentz-envelope}
b_{k-1}=\sup_{m\ge j}2^{-m/s}x_m\quad(2^j\le k<2^{j+1}),\qquad
2^{j/s}b_{2^j-1}\le\sum_{r\ge0}2^{-r/s}x_{j+r}.
\end{equation}
Then $b$ is nonincreasing and null. Minkowski's inequality for the summable geometric kernel in \eqref{eq:lorentz-envelope} gives $\|(2^{j/s}b_{2^j-1})_j\|_q\le(1-2^{-1/s})^{-1}\|x\|_q$, also for $q=\infty$. Hence $b\in\ell^{s,q}$ by \eqref{eq:lorentz-dyadic}. On the preceding block $2^{j-1}\le k<2^j$, $\mathcal A_k\ge\mathcal A_{2^j}$ and $b_{k-1}\ge2^{-j/s}x_j$, so
\begin{equation}\label{eq:lorentz-energy-lower}
\sum_nA_n^2b_n^2\ge\frac12\sum_{j\ge1}w_j(A)x_j^2.
\end{equation}
If $q\le2$ and $w$ is unbounded, choose $j_k\uparrow\infty$ with $w_{j_k}\ge2^{2k}$ and put $x_{j_k}=2^{-k}$, zero elsewhere. Then $x\in\ell^q$ but its weighted square sum diverges. For $2<q<\infty$, put $t=q/(q-2)$. If $\sum w_j^t=\infty$, choose disjoint finite intervals $I_k$ with $S_k:=\sum_{j\in I_k}w_j^t\ge2^{2kt}$ and define
\begin{equation}\label{eq:lorentz-dual-blocks}
x_j=S_k^{-1/2}w_j^{1/(q-2)}\ (j\in I_k),\quad
\sum_{j\in I_k}w_jx_j^2=1,\quad
\|x|_{I_k}\|_q=S_k^{-1/(2t)}\le2^{-k}.
\end{equation}
Set $x=0$ elsewhere. Formula~\eqref{eq:lorentz-dual-blocks} gives $x\in\ell^q$ with divergent energy. If $q=\infty$ and $\sum w_j=\infty$, take $x_j=1$. In all cases \eqref{eq:lorentz-energy-lower} supplies a witness against multiplier membership. The finite coordinate $j=0$ cannot affect these tests. This proves the first equivalence in \eqref{eq:lorentz-criterion}; Theorem~\ref{thm:solid}, using Lemma~\ref{lem:spaces}, proves the second and its additional assertions. The witnesses in \eqref{eq:lorentz-envelope} have positive coordinates whenever the energy diverges, so they also give direct admissible majorants $Ab$.
\end{proof}

Bounding $\mathcal A_k$ at both dyadic endpoints and shifting the lower bound by one block gives $\sum_jw_j(A)\asymp\sum_nA_n^2(n+1)^{-2/s}$. Thus $q=\infty$ recovers Theorem~\ref{thm:weak}. Compare the classical multiplier component of Theorem~\ref{thm:lorentz} with \cite{Carando2009}.

\section{Power, logarithmic, and endpoint profiles}\label{sec:examples}
Write $a_s=1/2-1/s$, so $a_2=0$ and $a_s=1/R_s$ for $s>2$. By Corollary~\ref{cor:comparison}, positive sequences $A_n\asymp B_n$ have the same diagonal images; the notation means two-sided comparison by constants; finite initial terms do not affect any convergence test below.

\begin{lemma}[Power--logarithmic summability]\label{lem:power-log}
For real $u,v$,
\begin{equation}\label{eq:power-log-test}
\sum_{n\ge2}n^{-u}(\log n)^{-v}<\infty
\quad\Longleftrightarrow\quad u>1\text{ or }[u=1\text{ and }v>1].
\end{equation}
\end{lemma}
\begin{proof}
For every $\eta>0$, $(\log n)^{|v|}\le n^\eta$ eventually, since $\log\log n/\log n\to0$. If $u>1$, choose $\eta<u-1$ and bound the terms above by $n^{-u+\eta}$. If $u<1$, choose $\eta<1-u$ and bound them below by $n^{-u-\eta}$. The power-series integral test gives the respective conclusions in \eqref{eq:power-log-test}. For $u=1$, the summand is eventually decreasing; substituting $t=\log x$ in its integral gives $\int t^{-v}dt$, finite exactly when $v>1$. This also verifies the equality boundary.
\end{proof}

\begin{corollary}[Monotone model profiles]\label{cor:models}
For positive nonincreasing $A$, the following statements hold.
\begin{enumerate}[label=\textup{(\alph*)}]
\item If $A_n\asymp(n+1)^{-\alpha}$, $\alpha\ge0$, then $D_A(\Es)\cap\U\ne\varnothing$ exactly when $\alpha\le a_s$. In particular, no positive power decay works at $s=2$.
\item If $A_n\asymp(\log(n+2))^{-\beta}$, $\beta>0$, or $A_n=\exp(-[\log(n+2)]^\theta)$, $0<\theta<1$, then $D_A(\Es)\cap\U\ne\varnothing$ for every $s\ge2$.
\item If $s>2$ and
\begin{equation}\label{eq:critical-profile}
A_n\asymp(n+1)^{-a_s}(\log(n+2))^{-\beta},\qquad \beta\in\mathbb R,
\end{equation}
then an $\Es$ factor works for every $\beta$, whereas a weak-space factor works exactly when $\beta\le1/2$. At $s=2$ the weak threshold holds for $A_n\asymp(\log(n+2))^{-\beta}$. For the subexponential profiles in \textup{(b)}, a weak factor is impossible at $s=2$ and exists at every $s>2$.
\end{enumerate}
\end{corollary}
\begin{proof}
Lemma~\ref{lem:power-log} gives $A\in\ell^r$ for powers exactly when $\alpha r>1$. Such $r<R_s$ exists exactly when $\alpha>a_s$, including $s=2$ by \eqref{eq:spaces}; Theorem~\ref{thm:monotone} proves (a). Logarithmic profiles belong to no finite $\ell^r$ by the same lemma. For the other profile in (b), eventually $r(\log(n+2))^\theta\le\log(n+2)$, so $A_n^r\ge(n+2)^{-1}$; \eqref{eq:monotone} proves (b). For \eqref{eq:critical-profile}, $a_sr<1$ whenever $r<R_s$, so every such power sum diverges. The weak energy is
\begin{equation}\label{eq:critical-energy}
A_n^2(n+1)^{-2/s}\asymp\frac1{(n+1)(\log(n+2))^{2\beta}}.
\end{equation}
By Lemma~\ref{lem:power-log}, \eqref{eq:critical-energy} diverges exactly when $2\beta\le1$; Theorem~\ref{thm:weak} proves (c), also for logarithmic models at $s=2$. For the subexponential profile at $s=2$, substitute $t=\log x$: its weak-energy integral becomes $\int e^{-2t^\theta}\,dt<\infty$, since eventually $e^{-2t^\theta}\le t^{-2}$. If $s>2$, choose $0<\varepsilon<1-2/s$; eventually $A_n^2\ge(n+2)^{-\varepsilon}$, giving divergent weak energy.
\end{proof}

\begin{corollary}[Critical Lorentz thresholds]\label{cor:lorentz-models}
Let $s>2$ and let $A$ be positive nonincreasing with \eqref{eq:critical-profile}. Then
\begin{equation}\label{eq:lorentz-log-threshold}
D_A(\ell^{s,q})\cap\U\ne\varnothing\iff
\begin{cases}\beta<0,&1\le q\le2,\\
\beta\le\tfrac12-\tfrac1q,&2<q<\infty,\\
\beta\le\tfrac12,&q=\infty.
\end{cases}
\end{equation}
For $s=2$ the same tests hold for the nonincreasing logarithmic models, necessarily with $\beta\ge0$.
\end{corollary}
\begin{proof}
Here $w_j(A)\asymp(j+1)^{-2\beta}$. The spaces $Y_q$ in \eqref{eq:lorentz-weight} therefore give \eqref{eq:lorentz-log-threshold} by the power test and Theorem~\ref{thm:lorentz}. To make the finite-$q$ equality case explicit, let $2<q<\infty$, $\beta=1/2-1/q$, and $1/q<\delta\le1/2$. On a sufficiently late tail put
\begin{equation}\label{eq:double-log-witness}
b_n=(n+1)^{-1/s}(\log(n+2))^{-1/q}(\log\log(n+3))^{-\delta}.
\end{equation}
Complete the initial segment constantly to obtain positive nonincreasing $b$. Substitution $u=\log\log x$ in the integral test gives
\begin{equation}\label{eq:double-log-tests}
\sum_n(n+1)^{q/s-1}b_n^q<\infty\ (\delta q>1),\qquad
\sum_nA_n^2b_n^2=\infty\ (2\delta\le1).
\end{equation}
Indeed the tail summands in \eqref{eq:double-log-tests} are respectively comparable to $[n\log n(\log\log n)^{\delta q}]^{-1}$ and $[n\log n(\log\log n)^{2\delta}]^{-1}$. Thus \eqref{eq:double-log-witness} is an explicit Lorentz factor witness, and Lemma~\ref{lem:pog} applied to $Ab$ realizes the equality case. The calculation is unchanged at $s=2$.
\end{proof}

\begin{remark}[Two endpoint distinctions]\label{rem:endpoints}
For $s>2$ take $A_n=(n+2)^{-1/R_s}(\log(n+2))^{-\beta}$ with $\beta>1/R_s$. Then
\begin{equation}\label{eq:endpoint-obstruction}
\sum A_n^{R_s}=\sum\frac1{(n+2)(\log(n+2))^{\beta R_s}}<\infty,
\qquad \|Ac\|_2\le\|A\|_{R_s}\|c\|_s.
\end{equation}
The estimate is \eqref{eq:holder} with $q=s$. Proposition~\ref{prop:l2} and \eqref{eq:endpoint-obstruction} exclude $\ell^s$ factors; they do not make $A$ a multiplier on the larger $\Es$. Corollary~\ref{cor:models}(c) supplies $\Es$ factors. Thus the open union in \eqref{eq:monotone} cannot be replaced by $\ell^{R_s}$. Also, \eqref{eq:weak-inclusion} is strict: a sequence eventually equal to $(n+2)^{-1/s}\log(n+2)$ belongs to every $\ell^q$, $q>s$, but its decreasing rearrangement has infinite weak quasi-norm. For the critical profiles with $\beta>1/2$, Corollary~\ref{cor:models} shows that every admissible $\Es$ factor lies outside $\ell^{s,\infty}$.
\end{remark}

\section{Obstructions for nonmonotone profiles}\label{sec:obstructions}
\begin{proposition}[Reflection obstruction]\label{prop:reflection}
The positive null profile
\begin{equation}\label{eq:reflection-profile}
A_0=1,\qquad A_{2k-1}=(\log(k+2))^{-1},\qquad A_{2k}=2^{-k}\quad(k\ge1)
\end{equation}
belongs to no finite $\ell^r$, but admits no mode of universality with any $c\in\Es$.
\end{proposition}
\begin{proof}
Both branches in \eqref{eq:reflection-profile} tend to zero; the logarithmic branch has divergent $r$-power sum for every finite $r>0$, by the argument in Corollary~\ref{cor:models}. Fix $c\in\Es$ and put $d=Ac$. Lemma~\ref{lem:spaces} gives boundedness of $c$, hence $\sum_k|d_{2k}|\le\|c\|_\infty\sum_k2^{-k}<\infty$. Reflection $R(x)=2\pi-x$ leaves constants and cosines unchanged and negates sines. Thus for every finite $P=\sum_{n\in I}e_nd_n\Phi_n$, $|e_n|\le1$,
\begin{equation}\label{eq:reflection-bound}
|P(x)-P(Rx)|\le\frac2{\sqrt\pi}\sum_{k\ge1}|d_{2k}|=:C_c<\infty.
\end{equation}
Every partial sum in any mode satisfies \eqref{eq:reflection-bound}. If it converges almost everywhere to finite $F$, reflection preserves null sets, so passing to the limit gives $|F(x)-F(Rx)|\le C_c$ almost everywhere. The target $F=x^{-1/2}\mathbf1_{(0,\pi)}$, zero elsewhere, has difference $x^{-1/2}$ on $(0,\pi)$, exceeding $C_c$ on a set of positive measure. It is finite almost everywhere and in $L^1$, so it excludes each target class under consideration.
\end{proof}

Proposition~\ref{prop:reflection} tests the full sine branch. Rotation averages give restrictions on every divisibility subsystem. Interpret functions on the circle and define, for $m\ge1$,
\begin{equation}\label{eq:arithmetic-projection}
\mathcal T_mF(x)=\frac1m\sum_{l=0}^{m-1}F(x+2\pi l/m),\quad
\mathcal RF(x)=F(-x),\quad \Pi_m^\pm=\tfrac12(I\pm\mathcal R)\mathcal T_m.
\end{equation}
The geometric sum $m^{-1}\sum_{l=0}^{m-1}e^{2\pi ikl/m}$ equals one for $m\mid k$ and zero otherwise. Consequently \eqref{eq:arithmetic-projection} retains the sine frequencies $mk$ for the minus sign and the cosine frequencies $mk$, together with the constant, for the plus sign.

\begin{proposition}[Projected targets]\label{prop:projected}
If $(d_{2mk})_{k\ge1}\in\ell^2$, every finite almost-everywhere target obtainable by subseries, signs, or rearrangements satisfies $\Pi_m^-F\in L^2$. If $(d_{2mk-1})_{k\ge1}\in\ell^2$, it satisfies $\Pi_m^+F\in L^2$. In the rearrangement case the projected limit is the fixed $L^2$ sum of the corresponding original subsystem.
\end{proposition}
\begin{proof}
If $P_N\to F$ almost everywhere, remove the finitely many preimages of its exceptional null set under the substitutions in \eqref{eq:arithmetic-projection}. Outside their union all arguments converge, so
\begin{equation}\label{eq:projected-limit}
\Pi_m^\pm P_N\longrightarrow\Pi_m^\pm F\quad\text{almost everywhere}.
\end{equation}
The respective square-summable subsystem converges in $L^2$ by \eqref{eq:orthogonal}, including the one constant term for the plus sign. The subsequence argument in Proposition~\ref{prop:l2} identifies its $L^2$ limit with \eqref{eq:projected-limit}. For permutations the same proof of unconditional convergence in Proposition~\ref{prop:l2} makes the limit independent of the order.
\end{proof}

\begin{theorem}[Arithmetic-subsystem obstruction]\label{thm:arithmetic}
For arbitrary positive $A$, if for some $m\ge1$ either sequence
\begin{equation}\label{eq:arithmetic-condition}
(A_{2mk-1})_{k\ge1},\qquad (A_{2mk})_{k\ge1}
\end{equation}
belongs to $\bigcup_{0<r<R_s}\ell^r$, no $c\in\Es$ yields any mode of universality.
\end{theorem}
\begin{proof}
Restriction preserves $\Es$ by Lemma~\ref{lem:spaces}. Theorem~\ref{thm:mult} applied to the relevant sequence in \eqref{eq:arithmetic-condition} makes its coefficients square-summable. Proposition~\ref{prop:projected} then applies. On $(-\pi/m,\pi/m)$ set
\begin{equation}\label{eq:arithmetic-target}
F_m^+(x)=|x|^{-1/2},\qquad F_m^-(x)=\operatorname{sgn}(x)|x|^{-1/2}\quad(x\ne0),
\end{equation}
and extend with period $2\pi/m$, assigning finite values at the finitely many singular points. These functions are finite almost everywhere, not in $L^2$, and their parity and period give $\Pi_m^\pm F_m^\pm=F_m^\pm$. Thus \eqref{eq:arithmetic-target} contradicts the respective projected restriction. The targets are finite, so the conclusion covers both rearrangement classes.
\end{proof}

\begin{corollary}[Two full branches do not suffice]\label{cor:two-branches}
The positive null profile
\begin{equation}\label{eq:two-branch-profile}
A_0=1,\qquad A_{2k-1}=A_{2k}=\begin{cases}2^{-k},&k\text{ even},\\(\log(k+2))^{-1},&k\text{ odd},\end{cases}
\end{equation}
has both full sine and cosine branches outside every finite $\ell^r$, but admits no mode with any $c\in\Es$.
\end{corollary}
\begin{proof}
The odd-frequency parts of \eqref{eq:two-branch-profile} have divergent logarithmic $r$-power sums. Its even-frequency parts are $2^{-2k}\in\ell^r$ for every $r>0$. Theorem~\ref{thm:arithmetic} with $m=2$ proves the assertion.
\end{proof}

\section{Independent regular branches and finite targets}\label{sec:branches}
Split the nonconstant indices into
\begin{equation}\label{eq:branches}
J_j=\{4k-4+j:k\ge1\},\qquad A_k^{(j)}=A_{4k-4+j}\quad(1\le j\le4),
\end{equation}
and assign the constant index to $J_3$ when partitioning $\mathbb N_0$. The order of these branches is odd-frequency cosines, odd-frequency sines, even-frequency cosines, and even-frequency sines. Put $\mathcal TF(x)=F(x+\pi)$ and use $\mathcal R$ from \eqref{eq:arithmetic-projection}. The commuting involutions give
\begin{equation}\label{eq:four-projections}
\Pi_j=\tfrac14(I+\varepsilon_j\mathcal R)(I+\tau_j\mathcal T),\qquad
((\varepsilon_j,\tau_j))_{j=1}^4=((1,-1),(-1,-1),(1,1),(-1,1)).
\end{equation}
Each $\Pi_j$ retains exactly $J_j$, and $\Pi_3$ also retains the constant. Direct expansion in \eqref{eq:four-projections} yields
\begin{equation}\label{eq:projection-resolution}
\Pi_i\Pi_j=0\ (i\ne j),\qquad\Pi_j^2=\Pi_j,\qquad\sum_{j=1}^4\Pi_j=I.
\end{equation}

\begin{lemma}[Four-branch obstruction]\label{lem:branch-obstruction}
If one branch $(d_{4k-4+j})_{k\ge1}$ is square-summable, the series has none of the universality properties.
\end{lemma}
\begin{proof}
Apply the proof of Proposition~\ref{prop:projected} to the finite substitutions in \eqref{eq:four-projections}. The projected finite target must be in $L^2$. For $0<\delta<\pi/4$, take $h=x^{-1/2}\mathbf1_{(0,\delta)}$ and $F=\Pi_jh$. By \eqref{eq:projection-resolution}, $\Pi_jF=F$. On $(0,\delta)$ the other three substituted arguments lie outside the support of $h$, so $F(x)=x^{-1/2}/4$. It is finite almost everywhere but not in $L^2$, a contradiction. The constant in $J_3$ contributes only a finite term.
\end{proof}

\begin{lemma}[Finite projection transfer]\label{lem:transfer}
Let $u$ be the Fourier coefficients of a real $g\in\bigcap_{1\le p<2}L^p$ and suppose $u\in\Uf$. For each $j$, the series with coefficients $u$ restricted to $J_j$ represents every target in $\Pi_j\mathcal S$ by each finite-target mode. Its coefficients are those of $\Pi_jg$, which belongs to all the same $L^p$ spaces.
\end{lemma}
\begin{proof}
Measure preservation and the triangle inequality give $\|\Pi_jg\|_p\le\|g\|_p$, $p\ge1$. Both involutions are self-adjoint under integration and commute. Thus, for each basis function in \eqref{eq:system},
\begin{equation}\label{eq:projection-coefficients}
\int(\Pi_jg)\Phi_n=\int g(\Pi_j\Phi_n)
=\begin{cases}u_n,&n\in J_j,\\0,&n\notin J_j,\end{cases}
\end{equation}
with the constant convention in \eqref{eq:branches}. This is a change of variables in four integrals, independent of Fourier-series convergence.

Choose a selected or sign-modified series for a finite target $F$. Projecting its partial sums preserves almost-everywhere convergence by the null-set argument of \eqref{eq:projected-limit}, and retains exactly the selected or signed $J_j$ terms. Its limit is $\Pi_jF$. For a rearrangement $\sigma$, let $k_1<k_2<\cdots$ be the positions at which $\sigma(k_m)\in J_j$, and put $\tau(m)=\sigma(k_m)$. The restriction retains all $J_j$ indices, including zero coefficients, and
\begin{equation}\label{eq:restricted-permutation}
\sum_{l=1}^m u_{\tau(l)}\Phi_{\tau(l)}
=\Pi_j\!\left(\sum_{l=0}^{k_m}u_{\sigma(l)}\Phi_{\sigma(l)}\right).
\end{equation}
The first position may be zero. Since $k_m\to\infty$, \eqref{eq:restricted-permutation} converges along all restricted partial sums. Together with \eqref{eq:projection-coefficients} this proves the lemma.
\end{proof}

\begin{theorem}[Independent four-branch criterion]\label{thm:branches}
Let $A_n>0$ and, for each $1\le j\le4$, suppose there is a positive nonincreasing model $B^{(j)}$ and constants $0<m_j\le M_j<\infty$ such that
\begin{equation}\label{eq:branch-comparison}
m_jB_k^{(j)}\le A_k^{(j)}\le M_jB_k^{(j)}\quad(k\ge1).
\end{equation}
No ordering between the models is assumed. For $1\le q\le\infty$, use $w$ and $Y_q$ from \eqref{eq:lorentz-weight}, with one-based model indexing $w_m(B)=2^{m(1-2/s)}B_{2^m}^2$. Then
\begin{align}
D_A(\Es)\cap\Uf\ne\varnothing
&\Longleftrightarrow B^{(j)}\notin\bigcup_{0<r<R_s}\ell^r\quad\text{for every }j,\label{eq:branch-Es}\\
D_A(\ell^{s,\infty})\cap\Uf\ne\varnothing
&\Longleftrightarrow\sum_{k\ge1}(B_k^{(j)})^2k^{-2/s}=\infty\quad\text{for every }j,\label{eq:branch-weak}\\
D_A(\ell^{s,q})\cap\Uf\ne\varnothing
&\Longleftrightarrow w(B^{(j)})\notin Y_q\quad\text{for every }j.\label{eq:branch-lorentz}
\end{align}
Each equivalence also holds if only one of the finite-target modes is required. All sufficient directions produce Fourier coefficients of one real function in every $L^p$, $1\le p<2$.
\end{theorem}
\begin{proof}
Restriction preserves $X$ by Lemma~\ref{lem:spaces}. If a branch model is a multiplier, \eqref{eq:branch-comparison} makes that branch of $Ac$ square-summable. Theorems~\ref{thm:mult}, \ref{thm:weak}, and~\ref{thm:lorentz}, followed by Lemma~\ref{lem:branch-obstruction}, prove necessity in \eqref{eq:branch-Es}--\eqref{eq:branch-lorentz}. For the weak case, writing $v_k=c_{4k-4+j}$,
\begin{equation}\label{eq:branch-energy-bound}
\sum_k|A_k^{(j)}v_k|^2\le M_j^2\|c\|_{s,\infty}^2\sum_k(B_k^{(j)})^2k^{-2/s},
\end{equation}
by Lemma~\ref{lem:rearrangement}; \eqref{eq:branch-energy-bound} verifies \eqref{eq:branch-weak} directly.

For sufficiency, repeat each model separately:
\begin{equation}\label{eq:repeated-models}
\widehat B_0^{(j)}=B_1^{(j)},\qquad\widehat B_n^{(j)}=B_{\lfloor(n+3)/4\rfloor}^{(j)}\quad(n\ge1).
\end{equation}
For $n=4k-4+j$, $1\le j\le4$, $\lfloor(n+3)/4\rfloor=k$, so $\widehat B_n^{(j)}=B_k^{(j)}$ exactly. For $r>0$, its $r$-power sum equals $(B_1^{(j)})^r+4\sum_k(B_k^{(j)})^r$. Also $4k-2\le n+1\le4k+1$ for $4k-3\le n\le4k$, so
\begin{equation}\label{eq:repeat-tests}
\sum_{n\ge0}(\widehat B_n^{(j)})^2(n+1)^{-2/s}\asymp\sum_{k\ge1}(B_k^{(j)})^2k^{-2/s}.
\end{equation}
For Lorentz spaces, $B\in\M(\ell^{s,q},\ell^2)$ if and only if its fourfold repetition is a multiplier: restrict and assemble the four coordinate classes in Lemma~\ref{lem:spaces}; the constant is harmless. Thus Theorem~\ref{thm:monotone}, \ref{thm:weak}, or~\ref{thm:lorentz}, according to $X$, applied using \eqref{eq:repeated-models}--\eqref{eq:repeat-tests} gives four real functions $f_j$ with strongly universal coefficients $d_n^{(j)}=\widehat B_n^{(j)}u_n^{(j)}$, $u^{(j)}\in X$.

Set $f=\sum_j\Pi_jf_j$. By Lemma~\ref{lem:transfer}, it is in all the stated $L^p$ spaces and its coefficients equal $d_n^{(j)}$ on $J_j$, with constant $d_0^{(3)}$. For $n=4k-4+j$, \eqref{eq:branch-comparison} and \eqref{eq:repeated-models} give
\begin{equation}\label{eq:branch-factor}
c_n=d_n^{(j)}/A_n,\qquad |c_n|\le m_j^{-1}|u_n^{(j)}|,
\qquad c_0=d_0^{(3)}/A_0.
\end{equation}
Lemma~\ref{lem:spaces} gives $c\in X$ from \eqref{eq:branch-factor}.

For finite $F$, Lemma~\ref{lem:transfer} supplies branch selectors or signs with limits $\Pi_jF$. At each natural upper index their four partial sums add; by \eqref{eq:projection-resolution} the limit is $F$. For permutations, interleave the four branch permutations from that lemma, taking one term from each cyclically. Every index occurs once and each local count tends to infinity. Thus at every intermediate position the sum converges to $\sum_j\Pi_jF=F$. Intersecting four full-measure sets for this target completes all equivalences on one Fourier coefficient sequence.
\end{proof}

\begin{remark}[The target-class boundary]\label{rem:finite-targets}
Theorem~\ref{thm:branches} uses finite targets so that $\Pi_jF$ and their sum are defined. This limitation is substantive for the method: on a four-point orbit, a character component infinite at one point is infinite at every point, with a possible sign change. A defined finite sum at another orbit point then requires all components finite. Consequently a target infinite at one point and finite at another cannot be decomposed into such extended-real components without an undefined cancellation. The measurable target $+\infty$ on $(0,\pi/4)$ and zero elsewhere supplies such orbits on a set of positive measure. This does not disprove extended-real universality by a different construction. Corollary~\ref{cor:extended} gives a sufficient additional model condition.
\end{remark}

\begin{corollary}[Finite residue-class decompositions]\label{cor:residues}
Fix $N\ge1$. Partition positive frequencies into $k\equiv\pm r\pmod N$, $0\le r\le\lfloor N/2\rfloor$, and split each class into cosine and sine branches. Suppose $A$ restricted to each branch, in increasing frequency order, is comparable to a positive nonincreasing model $B$. For $X=\Es$ or $\ell^{s,q}$, a Fourier coefficient sequence in $D_A(X)\cap\Uf$ exists if and only if every branch model is outside $\M(X,\ell^2)$. The latter tests are Theorems~\ref{thm:mult} and~\ref{thm:lorentz}. Each individual finite-target mode has the same criterion.
\end{corollary}
\begin{proof}
Let $T_\ell F(x)=F(x+2\pi\ell/N)$, $RF(x)=F(-x)$, and set
\begin{equation}\label{eq:residue-projections}
Q_r=\frac{v_r}{N}\sum_{\ell=0}^{N-1}\cos(2\pi r\ell/N)T_\ell,\quad
\Pi_{r,\pm}=\tfrac12(I\pm R)Q_r,\quad
v_r=\begin{cases}1,&2r\equiv0\pmod N,\\2,&\text{otherwise}.\end{cases}
\end{equation}
For clarity, on the exponential $e_k(x)=e^{ikx}$ one has
\[
 Q_re_k=e_k\,\frac{v_r}{2N}\sum_{\ell=0}^{N-1}
 \left(e^{2\pi i(r+k)\ell/N}+e^{2\pi i(k-r)\ell/N}\right).
\]
Each geometric sum is $N$ precisely when its exponent is divisible by $N$, and is zero otherwise. Thus $Q_r$ is the identity on the frequency pair $k\equiv\pm r\pmod N$ and zero on the other residue pairs. If $r=0$, or if $N$ is even and $r=N/2$, the two congruences coincide and the normalization $v_r=1$ removes the double counting; otherwise $v_r=2$. Since $RT_\ell R=T_{-\ell}$ and the cosine weights are unchanged by $\ell\mapsto-\ell$, $RQ_r=Q_rR$. The same finite geometric sums give $Q_rQ_t=0$ for $r\ne t$ and $\sum_rQ_r=I$. Thus \eqref{eq:residue-projections} gives complementary real projections, with the constant in $(0,+)$. These self-adjoint finite substitution sums are bounded on $L^p$, $p\ge1$; Lemma~\ref{lem:transfer}'s null-set and permutation arguments apply. For necessity, take $h=x^{-1/2}\mathbf1_{(0,\delta)}$, $0<\delta<\pi/N$. On $(0,\delta)$ only the identity substitution meets its support, so $\Pi_{r,\pm}h=(v_r/2N)x^{-1/2}\notin L^2$. This proves the analogue of Lemma~\ref{lem:branch-obstruction}.

For sufficiency list the nonconstant coefficient indices of one branch as $n_1<n_2<\cdots$. For a fixed residue/parity branch, admissible coefficient indices form a fixed finite union of residue classes modulo $4N$. Hence there are constants $c_N,C_N>0$ such that its counting function $N_J(t)=\#\{k:n_k\le t\}$ satisfies $c_Nt-C_N\le N_J(t)\le C_Nt+C_N$ for all sufficiently large $t$. Consequently $n_k\asymp k$. Moreover every interval of $4N$ consecutive coefficient indices contains a member of the branch, so $n_{k+1}-n_k\le4N$. Define $m(n)=\max(1,\#\{k:n_k\le n\})$ and $\widehat B_n=B_{m(n)}$. This model is nonincreasing, equals $B_k$ at $n_k$, and repeats each value at most $4N$ times apart from a harmless finite initial segment. Restriction and finite assembly in Lemma~\ref{lem:spaces} show that $B$ is a multiplier exactly when $\widehat B$ is. Theorem~\ref{thm:solid} realizes each lifted model. Project these realizations and interleave the finitely many branch permutations as in Theorem~\ref{thm:branches}. All partial sums have the required finite limits, completing the proof.
\end{proof}

\begin{remark}[Why a spectral idempotent alone is insufficient]\label{rem:idempotents}
The coset-ring description of idempotent measures (see \cite{GreenSanders2008}) does not replace the finite-substitution argument in Corollary~\ref{cor:residues}. Indeed, the Fej\'er polynomials $K_n(x)=n^{-1}[\sin(nx/2)/\sin(x/2)]^2$ tend to zero almost everywhere, while their normalized integrals equal one \cite[Vol.~I, Ch.~III]{Zygmund2002}. The projector onto the constant frequency therefore fails to preserve almost-everywhere convergence. No characterization of all admissible projection obstructions is claimed here.
\end{remark}

\section{Extended targets and ordered model families}\label{sec:extended}
\begin{corollary}[A weakest model and extended targets]\label{cor:extended}
Under \eqref{eq:branch-comparison}, suppose for some $j_0$ and $C<\infty$ that
\begin{equation}\label{eq:weakest}
B_k:=B_k^{(j_0)}\le CB_k^{(j)}\quad(k\ge1,\ 1\le j\le4).
\end{equation}
Then
\begin{align}
D_A(\Es)\cap\U\ne\varnothing
&\Longleftrightarrow B\notin\bigcup_{0<r<R_s}\ell^r,\label{eq:extended-Es}\\
D_A(\ell^{s,\infty})\cap\U\ne\varnothing
&\Longleftrightarrow\sum_{k\ge1}B_k^2k^{-2/s}=\infty,\label{eq:extended-weak}\\
D_A(\ell^{s,q})\cap\U\ne\varnothing
&\Longleftrightarrow w(B)\notin Y_q\quad(1\le q\le\infty).\label{eq:extended-lorentz}
\end{align}
Each condition characterizes existence for each mode separately, including extended-real rearrangements; the coefficients again come from one function in all the stated $L^p$ spaces.
\end{corollary}
\begin{proof}
Necessity follows from Theorem~\ref{thm:branches}, using inclusion~\eqref{eq:classes}. Define $\widehat B_0=B_1$, $\widehat B_n=B_{\lfloor(n+3)/4\rfloor}$ for $n\ge1$. This is nonincreasing, and \eqref{eq:branch-comparison}--\eqref{eq:weakest} give
\begin{equation}\label{eq:common-model}
A_n\ge a\widehat B_n\quad(n\ge0),\qquad
a=\min\{A_0/B_1,m_1/C,m_2/C,m_3/C,m_4/C\}>0.
\end{equation}
Here $C>0$ follows from positivity in \eqref{eq:weakest}. The repetition tests \eqref{eq:repeated-models}--\eqref{eq:repeat-tests}, bounded-repetition stability in Lemma~\ref{lem:spaces}, and Theorem~\ref{thm:monotone}, \ref{thm:weak}, or~\ref{thm:lorentz} give strongly universal coefficients $d_n=\widehat B_nu_n$, with $u$ in the corresponding space. By \eqref{eq:common-model}, $c_n=d_n/A_n$ satisfies $|c_n|\le a^{-1}|u_n|$. Lemma~\ref{lem:spaces} gives membership, proving \eqref{eq:extended-Es}, \eqref{eq:extended-weak}, and~\eqref{eq:extended-lorentz}. Necessity already holds for a single finite-target mode by Theorem~\ref{thm:branches}; sufficiency supplies all modes with the extended class, completing the additional statements.
\end{proof}

\begin{corollary}[Four power--logarithmic branches]\label{cor:four-models}
Suppose $A_n>0$ and
\begin{equation}\label{eq:four-power-log}
A_k^{(j)}\asymp(k+1)^{-\alpha_j}(\log(k+2))^{-\beta_j},\qquad
\alpha_j\ge0,\quad\beta_j\in\mathbb R,
\end{equation}
where $\beta_j\ge0$ if $\alpha_j=0$. Put
\begin{equation}\label{eq:extreme-parameters}
\alpha_* =\max_j\alpha_j,\qquad
\beta_* =\max\{\beta_j:\alpha_j=\alpha_*\},\qquad a_s=\tfrac12-\tfrac1s.
\end{equation}
Then
\begin{align}
D_A(\Es)\cap\U\ne\varnothing&\Longleftrightarrow\alpha_*\le a_s,\label{eq:power-Es}\\
D_A(\ell^{s,\infty})\cap\U\ne\varnothing
&\Longleftrightarrow[\alpha_*<a_s]\text{ or }[\alpha_*=a_s\text{ and }\beta_*\le\tfrac12].\label{eq:power-weak}
\end{align}
For $\ell^{s,q}$ the condition is $\alpha_*<a_s$, or $\alpha_*=a_s$ with the corresponding $\beta_*$ test in \eqref{eq:lorentz-log-threshold}. The same criteria apply to each mode separately.
\end{corollary}
\begin{proof}
For positive $\alpha_j$, the logarithmic derivative of the model in \eqref{eq:four-power-log} is eventually negative; for $\alpha_j=0$ it is nonpositive by $\beta_j\ge0$. Replace finite initial segments to get comparable positive nonincreasing models. The one with parameters \eqref{eq:extreme-parameters} is dominated by all the others: a larger power gives a smaller sequence up to constants, regardless of logarithms, and at equal power the larger logarithmic exponent gives the smaller sequence. Thus \eqref{eq:weakest} holds.

For $s>2$ the deciding model belongs to some $\ell^r$, $r<R_s$, precisely when $\alpha_*>a_s$, by the power test in Corollary~\ref{cor:models}. At equality every such power exponent is less than one. For $s=2$ the condition is $\alpha_*>0=a_2$. This proves \eqref{eq:power-Es} from Corollary~\ref{cor:extended}. Its weak energy has the form
\begin{equation}\label{eq:four-model-energy}
\sum_{k\ge1}k^{-2\alpha_*-2/s}(\log(k+2))^{-2\beta_*}.
\end{equation}
Lemma~\ref{lem:power-log} shows that \eqref{eq:four-model-energy} diverges exactly in the cases in \eqref{eq:power-weak}. Corollary~\ref{cor:extended} gives that equivalence. For Lorentz spaces, \eqref{eq:lorentz-weight} applied to the weakest model has $w_m\asymp2^{2m(a_s-\alpha_*)}(m+1)^{-2\beta_*}$. It tends exponentially to infinity if $\alpha_*<a_s$, belongs to every $Y_q$ if $\alpha_*>a_s$, and at equality is decided by Corollary~\ref{cor:lorentz-models}. This proves all remaining statements.
\end{proof}

\begin{remark}[Interlacing jumps and failed global tests]\label{rem:branch-examples}
For $s=4$, $a_s=1/4$, define the branches by the exact models \eqref{eq:four-power-log}. The following three choices illustrate Corollary~\ref{cor:four-models}:
\begin{equation}\label{eq:branch-examples}
\begin{aligned}
\textup{(i)}\quad &\alpha=(1/8,1/4,1/8,1/4),&&\beta=(0,0,0,0);\\
\textup{(ii)}\quad &\alpha=(1/8,1/8,1/8,3/8),&&\beta=(0,0,0,0);\\
\textup{(iii)}\quad &\alpha=(1/4,1/4,1/4,1/4),&&\beta=(0,0,0,1).
\end{aligned}
\end{equation}
In \eqref{eq:branch-examples}(i) both spaces admit strong universality. Yet $A_{4k-1}/A_{4k-2}=(k+1)^{1/8}\to\infty$. Comparability $mD_n\le A_n\le MD_n$ to any nonincreasing $D$ would bound each consecutive increase by $M/m$, so no such comparison is possible. In (ii) three slower branches make $A$ outside every $\ell^r$, $r<4$, while the fourth branch has $\alpha_*>1/4$ and forbids universality. In (iii) an $\Es$ factor works, but no weak-space factor does; the global energy $\sum A_n^2(n+1)^{-1/2}$ still diverges because of the first three branches. Thus neither the global multiplier test nor the global weak energy decides nonmonotone universality. Unlike Corollary~\ref{cor:two-branches}, these profiles have regular branches, so their exact behavior is determined by the branch criterion.
\end{remark}

\section{Realization beyond a single monotone majorant}\label{sec:separation}
Theorem~\ref{thm:branches} can apply when every compatible monotone majorant is square-summable.

\begin{proposition}[Compatible majorants with bounded relative cost]\label{prop:envelope}
For any positive profile $A$, define its nonincreasing lower envelope by
\begin{equation}\label{eq:envelope}
A_n^{\downarrow}=\min_{0\le l\le n}A_l.
\end{equation}
There exists a positive nonincreasing null sequence $\rho\notin\ell^2$ with $\sup_n\rho_n/A_n<\infty$ if and only if $A^{\downarrow}\notin\ell^2$. Every nonnegative nonincreasing $\rho$ with $\rho\le CA$ satisfies
\begin{equation}\label{eq:envelope-domination}
\rho_n\le CA_n^{\downarrow}\quad(n\ge0).
\end{equation}
\end{proposition}
\begin{proof}
For $l\le n$, monotonicity gives $\rho_n\le\rho_l\le CA_l$. Taking the minimum proves \eqref{eq:envelope-domination}, and hence necessity by \eqref{eq:envelope}. Conversely suppose $A^{\downarrow}\notin\ell^2$. Its positive nonincreasing coordinates have a limit $L\ge0$. If $L=0$, take $\rho=A^{\downarrow}$; it has all the required properties and $\rho/A\le1$. If $L>0$, take $\rho_n=(n+1)^{-1/2}$, whose square sum diverges and whose ratio to $A$ is bounded by $1/L$. These two cases prove sufficiency. The equivalence concerns bounded ratios; it does not claim membership of the ratio in $\Es$.
\end{proof}

\begin{theorem}[Separation from compatible monotone majorants]\label{thm:separation}
There is a positive null profile $A$ satisfying the regular-branch assumptions \eqref{eq:branch-comparison} such that
\begin{equation}\label{eq:separation-realization}
D_A(\Es)\cap\Uf\ne\varnothing\quad\text{for every }2\le s<\infty,
\end{equation}
In fact one fixed $c\in E_2\subset\Es$ gives \eqref{eq:separation-realization} for all $s\ge2$. Every nonnegative nonincreasing $\rho$ with $\sup_n\rho_n/A_n<\infty$ belongs to $\ell^1$. In particular, no positive nonincreasing $\rho\notin\ell^2$ can satisfy $\rho/A\in\Es$, and no branch model is dominated by all the other models up to a fixed constant.
\end{theorem}
\begin{proof}
We first construct two positive nonincreasing null sequences $B^{(1)},B^{(2)}$ indexed by $k\ge1$. Set initial heights $b_{1,0}=b_{2,0}=1$ and $K_1=1$. On stage $j\ge1$, let $\iota_j=1$ for odd $j$ and $\iota_j=2$ for even $j$, and put
\begin{equation}\label{eq:separation-recursion}
\begin{gathered}
h_j=b_{\iota_j,j-1},\quad L_j=\lceil h_j^{-j}\rceil,\quad K_{j+1}=K_j+L_j,\\
t_j=\min\{b_{3-\iota_j,j-1}/2,\ h_j/2,\ 2^{-j}/L_j\},\quad
b_{\iota_j,j}=h_j,\quad b_{3-\iota_j,j}=t_j.
\end{gathered}
\end{equation}
Define $B_k^{(i)}=b_{i,j}$ on $K_j\le k<K_{j+1}$. The positive finite recursion \eqref{eq:separation-recursion} leaves each sequence unchanged on active stages and reduces it by at least half on every other stage. Both are therefore nonincreasing and null. Also $0<h_j\le1$ and $t_j\le h_j/2$, so the smaller height on stage $j$ is $t_j$.

For fixed $r>0$, each sequence has infinitely many active stages with $j>r$. Its contribution on such a stage satisfies
\begin{equation}\label{eq:separation-sums}
L_jh_j^r\ge h_j^{r-j}\ge1,\qquad
\sum_{k\ge1}\min(B_k^{(1)},B_k^{(2)})=\sum_{j\ge1}L_jt_j\le\sum_{j\ge1}2^{-j}<\infty.
\end{equation}
Consequently both $B^{(i)}$ lie outside every finite $\ell^r$, while their pointwise minimum belongs to $\ell^1$. 

Set
\begin{equation}\label{eq:separation-profile}
A_0=1,\qquad A_{4k-3}=A_{4k-1}=B_k^{(1)},\qquad
A_{4k-2}=A_{4k}=B_k^{(2)}\quad(k\ge1).
\end{equation}
This profile is positive and null; its four exact branch models are $B^{(1)},B^{(2)},B^{(1)},B^{(2)}$. The first assertion of \eqref{eq:separation-sums} verifies every condition in \eqref{eq:branch-Es}, for any $s\ge2$, so Theorem~\ref{thm:branches} at $s=2$ gives one $c\in E_2\subset\Es$ proving \eqref{eq:separation-realization} simultaneously for all $s\ge2$.

In the preceding four-element block choose an index where $A$ equals the smaller of its two model values. Formula~\eqref{eq:envelope} then gives $A_n^{\downarrow}\le\min(B_{k-1}^{(1)},B_{k-1}^{(2)})$ for $k\ge2$, $4k-3\le n\le4k$. For every nonnegative nonincreasing $\rho\le CA$, Proposition~\ref{prop:envelope} consequently gives
\begin{equation}\label{eq:separation-majorants}
\rho_n\le C\min(B_{k-1}^{(1)},B_{k-1}^{(2)}),\qquad
\sum_{n\ge5}\rho_n\le4C\sum_{k\ge1}\min(B_k^{(1)},B_k^{(2)})<\infty.
\end{equation}
The remaining initial terms are finite, so \eqref{eq:separation-majorants} proves $\rho\in\ell^1\subset\ell^2$, using Lemma~\ref{lem:spaces}. A ratio $\rho/A\in\Es$ is bounded by the same lemma and hence falls under this conclusion. Finally, if $B^{(1)}\le CB^{(2)}$, their minimum is at least $B^{(1)}/\max(1,C)$, contradicting its summability and $B^{(1)}\notin\ell^1$ from \eqref{eq:separation-sums}. The reverse comparison is equally impossible. Since the four models in \eqref{eq:separation-profile} repeat these two, no weakest model exists.
\end{proof}

\begin{remark}\label{rem:scope-separation}
Theorem~\ref{thm:separation} concerns $\Uf$; extended-real rearrangement targets remain open. Its obstruction covers every compatible monotone majorant, while Theorem~\ref{thm:branches} projects four separate realizations before assembly, without requiring such a global majorant.
\end{remark}

\section{A Walsh--Paley analogue}\label{sec:walsh}
On $[0,1)$ use Walsh--Paley functions $W_n=\prod_{j\ge0}r_j^{n_j}$, where $n=\sum_jn_j2^j$ and $r_j(x)=(-1)^{\lfloor2^{j+1}x\rfloor}$. A fixed $d$ is ternary-universal if every finite almost-everywhere measurable target is the almost-everywhere sum of $\sum\eta_nd_nW_n$ in natural order, with $\eta_n\in\{0,\pm1\}$. By \cite[Theorem~1, pp.~37--38]{GevorkyanNavasardyan1999}, every positive decreasing null $d\notin\ell^2$ has this property. We add Fourier origin and factor control, without claiming separate selector, sign, or permutation modes. Related Walsh results concern limit functions \cite{EpiskoposianMuller2014}, universal modification pairs \cite{EpiskoposianPair2024}, $L^p$-greedy universality after correction on small sets \cite{EpiskoposyanGreedy2022}, and weighted representation \cite{GrigorianEpiskoposian2001}. Here the target is represented almost everywhere with a prescribed diagonal cost. Fourier origin under the coefficient decay below also follows from the classical monotone-coefficient criterion recalled in \cite[p.~38]{GevorkyanNavasardyan1999}; we give a direct proof.

\begin{lemma}[Truncation and Fourier origin]\label{lem:walsh-truncation}
If $a_n>0$ is nonincreasing and $\sum a_n^2=\infty$, there is a strictly decreasing positive $d$ with $d_n\le2a_n$, $d\notin\ell^2$, and $d_n=O((n+1)^{-1/2})$. Its Walsh--Paley series converges almost everywhere and in every $L^p$, $1\le p<2$, to one function having Fourier--Walsh coefficients $d$.
\end{lemma}
\begin{proof}
For $n\ge1$ put $t_n=\min(a_n,n^{-1/2})$ and $t_0=t_1$. If $a_n\ge n^{-1/2}$ only finitely often, $t=a$ on a tail. Otherwise choose such indices $n_i\to\infty$ with $n_{i+1}>2n_i$. On each disjoint interval $\lceil n_i/2\rceil\le m\le n_i$, both $a_m$ and $m^{-1/2}$ are at least $n_i^{-1/2}$, so its squared $t$-sum is at least $1/2$. Thus $\sum t_n^2=\infty$ in either case. Taking $d_n=(1+2^{-n})t_n$ proves strict decrease, $d_n\le2a_n$, and all sequence assertions.

For completeness we prove Fourier origin directly. The Walsh Dirichlet kernels $D_N=\sum_{n<N}W_n$ satisfy \cite[p.~39]{GevorkyanNavasardyan1999}
\begin{equation}\label{eq:walsh-kernels}
D_{2^j}=2^j\mathbf1_{[0,2^{-j})},\qquad
D_{2^j+m}=D_{2^j}+W_{2^j}D_m\quad(0\le m<2^j).
\end{equation}
The first identity also follows from $D_{2^j}=\prod_{l<j}(1+r_l)$; the second follows from the binary definition. We record the Dirichlet-kernel estimate explicitly. Write
\[
 N=2^{j_1}+\cdots+2^{j_r},\qquad j_1>\cdots>j_r\ge0.
\]
Repeated use of the second identity in \eqref{eq:walsh-kernels} gives
\[
 D_N(x)=\sum_{\nu=1}^r \varepsilon_\nu(x)D_{2^{j_\nu}}(x),\qquad |\varepsilon_\nu(x)|=1.
\]
Hence, outside the dyadic endpoints,
\begin{equation}\label{eq:walsh-dirichlet-bound}
 |D_N(x)|\le \sum_{\nu=1}^r2^{j_\nu}\mathbf 1_{[0,2^{-j_\nu})}(x)
 \le \sum_{j:\,x<2^{-j}}2^j\le \frac{2}{x}.
\end{equation}
For $L\ge M$, Abel summation, written for the shifted tail, gives
\[
 \sum_{n=M}^{L}d_nW_n
 =d_L(D_{L+1}-D_M)+\sum_{n=M}^{L-1}(d_n-d_{n+1})(D_{n+1}-D_M).
\]
Using \eqref{eq:walsh-dirichlet-bound}, $d_n\downarrow0$, and telescoping of the coefficient differences, we obtain
\[
 \left|\sum_{n=M}^{L}d_nW_n(x)\right|\le \frac{8d_M}{x}\longrightarrow0
 \qquad(M\to\infty)
\]
for every fixed non-dyadic $x\in(0,1)$. Thus the Walsh series converges pointwise almost everywhere. If $d_n\le C(n+1)^{-1/2}$, split a partial sum at $M=\lceil1/x\rceil$. Its initial part is at most $2C\sqrt M$, while the remaining tail is at most $8d_M/x$. Since $0<x<1$ gives $x^{-1}\le M< x^{-1}+1\le2x^{-1}$, we have $2C\sqrt M\le2\sqrt2\,Cx^{-1/2}$. Also $M\ge x^{-1}$ implies $d_M\le C(M+1)^{-1/2}\le Cx^{1/2}$, and hence $8d_M/x\le8Cx^{-1/2}$. Consequently
\begin{equation}\label{eq:walsh-domination}
\sup_N\left|\sum_{n<N}d_nW_n(x)\right|\le(2\sqrt2+8)Cx^{-1/2}\quad(0<x<1\text{ a.e.}).
\end{equation}
Since $x^{-p/2}$ is integrable for $p<2$, dominated convergence in \eqref{eq:walsh-domination} gives $L^p$ convergence to the same pointwise limit $f$. In particular, $L^1$ convergence and $|W_k|=1$ permit passage under $\int\!\cdot\,W_k$; orthonormality gives $\int_0^1fW_k=d_k$ for every $k$.
\end{proof}

\begin{theorem}[Walsh diagonal criterion]\label{thm:walsh}
Let $X\subset c_0$ be a solid linear space stable under decreasing rearrangement, and let $A_n>0$ be nonincreasing. There exists $c\in X$ such that $d=Ac$ is ternary-universal and consists of the Fourier--Walsh coefficients of one real $f\in\bigcap_{1\le p<2}L^p[0,1)$ if and only if $A\notin\M(X,\ell^2)$. In particular the exact conditions are
\begin{equation}\label{eq:walsh-criteria}
X=\Es:\ A\notin\bigcup_{0<r<R_s}\ell^r;\qquad
X=\ell^{s,q}:\ w(A)\notin Y_q\quad(1\le q\le\infty).
\end{equation}
The chosen $d$ can be positive and strictly decreasing.
\end{theorem}
\begin{proof}
If $A$ is a multiplier, every ternary modification has square-summable coefficients. The orthonormality argument of Proposition~\ref{prop:l2}, now for $W_n$ and the target $x^{-1/2}$ on $(0,1)$, excludes universality. Otherwise choose $v\in X$ with $Av\notin\ell^2$ and set $b=v^*$. By Lemma~\ref{lem:rearrangement}, $Ab\notin\ell^2$; thus $b$ is positive, nonincreasing, null, and in $X$. Apply Lemma~\ref{lem:walsh-truncation} to $a=Ab$. Its sequence $d$ is ternary-universal by \cite[Theorem~1]{GevorkyanNavasardyan1999}, and has the asserted Fourier origin. Since $|d_n/A_n|\le2b_n$, solidity gives $c=d/A\in X$. Theorems~\ref{thm:mult} and~\ref{thm:lorentz} identify the conditions in \eqref{eq:walsh-criteria}. This proof uses no trigonometric realization theorem.
\end{proof}

\section{Scope and remaining questions}
The criteria, endpoints, and separation are Theorems~\ref{thm:lorentz}, \ref{thm:branches}, \ref{thm:walsh}, \ref{thm:separation} and Corollary~\ref{cor:lorentz-models}; their scope is clarified in Remarks~\ref{rem:HY}, \ref{rem:endpoints}, \ref{rem:branch-examples}, \ref{rem:scope-separation}. Residue projections (Corollary~\ref{cor:residues}, Remark~\ref{rem:idempotents}) do not resolve independent-branch extended targets beyond Corollary~\ref{cor:extended} (Remark~\ref{rem:finite-targets}), or arbitrary profiles subject to Theorem~\ref{thm:arithmetic}.

\paragraph{Declarations}
No specific funding was received. The author declares no competing interests. No research data were used.
\paragraph{Declaration of generative AI and AI-assisted technologies}
ChatGPT (OpenAI) was used for language editing, formatting assistance, exploratory drafting, and identifying points for additional author checking. All mathematical statements, proofs, references, revisions, and final decisions were reviewed and accepted by the author, who takes full responsibility for the manuscript.

\end{document}